\documentclass[authoryear,3p, times,10pt,final]{elsarticle}%gulliver
\usepackage{booktabs}
\usepackage{amsmath}
\usepackage{hyperref}
\usepackage{soul}

\usepackage{amsmath}
\begin{document}

\begin{frontmatter}

\newtheorem{thm}{Theorem}
\newtheorem{lem}[thm]{Lemma}
\newdefinition{rmk}{Remark}
\newdefinition{note}{Note}
\newproof{pf}{Proof}
\newproof{pot}{Proof of Theorem \ref{thm2}}

%%\usepackage{draftwatermark}
%\usepackage{draftwatermark}
%\SetWatermarkText{Confidential}
%%\SetWatermarkColor{red}
%\usepackage{ulem}
%%\allowdisplaybreaks
%\usepackage{xcolor}
%\usepackage[natbibapa]{apacite}
%\usepackage{comment}
%----------------------------------------------------------------------------------------
%	ARTICLE INFORMATION
%----------------------------------------------------------------------------------------
%\author[LIO]{Omar Kemmar}
%\author[LIO]{Karim Bouamrane}
\author[ITBA]{Rahimeh N Monemi \corref{cor1}}

\cortext[cor1]{Corresponding author, r.n.monemi@gmail.com}
%\address[SOTON]{Centre for Operational Research,  Management Science and Information Systems (CORMSIS), University of Southampton, Southampton, SO17 1BJ, United Kingdom}
%\address[OSM]{ Portsmouth Business School, Portsmouth, United Kingdom}

% %\address[LGI2A]{Laboratoire de G\'enie Informatique et d'Automatique, Universit\'{e} Artois, F-62400, B\'{e}thune, France}

%\address[LGI2A]{LGI2A (EA 3926), Universit\'{e} d'Artois, F-62400 B\'{e}thune, France}
%\address[LIO]{Laboratoire d'informatique d'Oran, Université Oran1, BP 1524 EL Mnaouer Oran, Algerie}
\address[ITBA]{IT and Business Analytics,  London, UK}

 %\address[LGI2A]{Laboratoire de G\'enie Informatique et d'Automatique, Universit\'{e} Artois, F-62400, B\'{e}thune, France}
% \address[NORD]{Univ Lille Nord de France, F-59000 Lille, France}
 %\address[ISIMA]{ISIMA, Universit\'{e} Blaise-Pascal, BP 10125, F-63173 Aubi\`{e}re Cedex, France}
 %\address[UBP]{LIMOS, UMR 6158-CNRS Universit\`{e} Blaise-Pascal, BP 10125, F-63173 Aubi\`{e}re Cedex, France}
% \address[IFFSTAR]{Ifsttar, Univ. Lille Nord de France, rue E\'lise\'e Reclus 20, 59666 Villeneuve d'Ascq, France}
% \address[CRIStAL]{CRIStAL, UMR 9189-CNRS, Ecole Centrale de Lille, 59651 Villeneuve d'Ascq, France}

 %\address[ITWM]{Fraunhofer Institute for Industrial Mathematics (ITWM), Kaiserslautern, Germany}

\title{A note on 'A multi-compartment VRP model for the health care waste transportation
problem'}
%----------------------------------------------------------------------------------------
%	ABSTRACT
%----------------------------------------------------------------------------------------

\begin{abstract}
A mathematical model and a genetic algorithm, referred to as an adaptive one, have been proposed in the paper by Nasreddine Ouertani, Hajer Ben-Romdhane, Issam Nouaouri, Hamid Allaoui, and Saoussen Krichen, titled 'A multi-compartment VRP model for the healthcare waste transportation problem,' published in the Journal of Computational Science in 2023 (72), pages 102-104. This model addresses the problem of waste disposal and the transportation of waste from healthcare facilities to treatment centers. In this note, we demonstrate that the model contains several minor and major flaws in its structure, making it, in short, incorrect. Therefore, we conclude that it cannot serve as a reliable benchmark for evaluating the proposed heuristic. We recommend amending the model to rectify some of these flaws.\end{abstract}
\end{frontmatter}
%----------------------------------------------------------------------------------------

%----------------------------------------------------------------------------------------
%	ARTICLE CONTENTS
%----------------------------------------------------------------------------------------
%
\section*{Introduction} % The \section*{} command stops section numbering

The transportation of medical waste generated by healthcare facilities (HCFs) to dedicated disposal centers is frequently conceptualized as a variation of the capacitated vehicle routing problem (CVRP). The CVRP itself is a well-studied problem, with a substantial body of literature available. However, due to the nature of medical waste, which comprises various types that cannot be mixed, it becomes apparent that the multi-compartment vehicle routing problem literature is applicable in modeling this issue, as observed in \cite{OUERTANI2023102104}.\\

In the following sections, we will first provide a problem description, the model, and its components as detailed in \cite{OUERTANI2023102104}. Afterward, we will conduct a thorough analysis of the existing shortcomings, and in the end, we will propose potential remedies.

\section{Problem Description}
Let $G=(N,A)$ be an undirected complete graph, where $N= N_c \cup {0}$ is the set of graph nodes, $N_c=\{1,\dots,N\}$  is the set of healthcare facilities (HCFs), while node 0 represents the depot. The set of edges $E = \{(i,j): i,j \in N, i\neq j\}$ denotes the distance; thus, a travel distance $d_{ij}, i.j\in N_c$ between two HCFs is associated to each edge $(i,j)\in E$. The collection and transportation tasks are ensured via a fleet of $K$ homogeneous multi-compartment vehicles. Each vehicle comprises $|P|$  compartments, each with a capacity of $Q_p(p=1,\dots, |P|)$. In other words, all vehicles contain $|P|$ compartments, where $p\leq |P|$, with P being equal to the number of types of waste to be collected. All vehicles have the same number of compartments.
Moreover, the fleet of vehicles is located at the collection center (depot) and must be returned to it after finishing the routing plan with the  minimum travel distance. Each HCF $i$, where $i\in N_c$  must be visited exactly once by only one vehicle. Moreover, a non-negative quantity of waste of type $p$ is designated to $q_{ip}$ to be taken from HCF $i$, where $q_{ip} \leq Q_p$.

\section{Modeling and Problem Settings}

The following parameters and variables have been used in the model in \cite{OUERTANI2023102104}:
\begin{table}[h]
\caption{Parameters.}
\label{tab:params}
\centering
\scalebox{0.8}{
\begin{tabular}{l l}
\toprule
%%\tabhead{Por 1} & \tabhead{Por 3}    \\
%%\midrule
${N}$:& Set of HCFs including the collection center,\\
${N_c}$:& Set of HCFs,\\
$K$: & Set of vehicles, \\
$P$: & Set of compartments,\\
$q_{ip}$:& Quantity of waste type $p$, where $p=1,\dots,|P|$ to be picked up from HCF  $i=q,\dots, N_c$,\\
$d_{ij}$: &Travel distance between HCF$_i$, HCF$_j$,\\
$Q_p$: & Capacity of compartment reserved to product, $p=1,\dots,|P|$\\
$g_j$: &  penalty cost for trailer $j$ that is not served at all.\\
\bottomrule\\
\end{tabular}
}
\end{table}

\begin{table}[h]
\caption{Decision variables.}
\label{tab:vars}
\centering
\scalebox{0.8}{
\begin{tabular}{l l}
\toprule
%%\tabhead{Por 1} & \tabhead{Por 3}    \\
%%\midrule
$x_{ijk}$:& 1 if vehicle $k$ visits HCF $j\in N_c$ after visiting HCF  $i\in N_c$, 0 otherwise,\\
$z_{ikp}$:& 1 if waste of type $p$  is loaded in vehicle $k$ from HCF $i\in N_c$, 0 otherwise,\\
$y_{ik}$: & 1 if vehicle $k$ visits HCF $i\in N_c$, 0 otherwise,\\
$u_{ikp}$:&  an integer variable representing the total carried quantity of \\
&   waste $p$   in vehicle $k$ after leaving HCF $i\in N_c$\\
\bottomrule\\
\end{tabular}
}
\end{table}

The mathematical model in \cite{OUERTANI2023102104} follows:
\begin{align}
\min & \sum_{k\in K}\sum_{i\in N}\sum_{j\in N} d_{ij}x_{ijk} \label{obj}\\
&\sum_{k\in K} y_{ik} = 1 & \forall i\in N_c \label{eq1}\\
&\sum_{k\in K} y_{0k} \leq |K| \label{eq2}\\
& \sum_{i\in N, i\neq j}x_{ijk} = y_{jk}, & \forall j\in N_c, k\in K \label{eq3}\\
& \sum_{j\in N, i\neq j}x_{ijk} = y_{ik}, & \forall i\in N_c, k\in K \label{eq4}\\
& u_{ikp}-u_{jkp} + Q_p x_{ijk} \leq Q_p - q_{jp} & \forall i,j \in N_c, p\in P,\ k\in K, i\neq j \label{eq5}\\
&q_{ip} \leq u_{ikp} \leq Q_p  & \forall i\in N_c, k\in K, p\in P\label{eq6}\\
&z_{jkp} \leq \sum_{i\in N, j\neq i} x_{ijk} &  \forall j\in N_c, k\in K, p\in P \label{eq7}\\
&\sum_{k\in K} z_{jkp} = 1 & \forall j \in N_c, p\in P  \label{eq8}\\
&\sum_{j\in N, j\neq 0} z_{jkp} q_{jp} \leq Q_p &  \forall k\in K, p\in P \label{eq9}\\
&x_{ijk}, z_{jkp}, y_{ik} \in \{0,1\} & \forall i,j\in N, p\in P, k\in K \label{eq10}
  \end{align}

The authors in \cite{OUERTANI2023102104} describe the model as follows: The objective function is a minimizes the total travel distance. Constraints \eqref{eq1}-\eqref{eq2} ensure that exactly one HCF is visited once by one vehicle and that all vehicles start from the collection center. Constraints \eqref{eq3}-\eqref{eq4} ensure that a vehicle serving HCF $i$ must arrive and leave from it. Constraints \eqref{eq5}-\eqref{eq6} represent the sub-tour elimination and impose, respectively, the requirement of capacity and connectivity between two HCFs. Constraint \eqref{eq7} set $z_{jkp}$ to zero for each HCF type $p$ if HCF $j$ is not visited by vehicle $k$. Constraint \eqref{eq8} ensure that each HCF of type $p$ is collected by exactly one vehicle. Constraint \eqref{eq9} check the compartment capacity constraint of each vehicle. Constraint \eqref{eq10} indicate the binary nature of the decision variables.

\section{The Flaws in the model \cite{OUERTANI2023102104} }

Let $N_c = {1, \dots, 10}$ and $N = N_c \cup {0}$. Also, let $K = {1, 2}$ and $P = {1, 2, 3}$. Given that the fleet is homogeneous, we assume that the capacity of each compartment for every waste type is the same and equals 500. This large capacity is intentionally set to avoid any infeasibility.\\

Additionally, we assume that the $q_{ip}$ values are given by the following matrix:
\begin{equation*}
\footnotesize
\begin{bmatrix}
7 & 7 & 6\\
8 & 5 & 4\\
4 & 6 & 3\\
3 & 2 & 2\\
6 & 4 & 3\\
6 & 1 & 3\\
8 & 7 & 7\\
9 & 9 & 8\\
0 & 4 & 7\\
2 & 3 & 2\\
\end{bmatrix}
\end{equation*}

and the distance table is given by:
\begin{equation*}
\footnotesize
\left [
  \begin{array}{ccccccccccc}
0 & 10 & 18 & 11 & 10 & 10 & 11 & 17 & 14 & 12 & 19 \\
11 &  0&17 & 12 & 15 & 18 & 11 & 16 & 14 & 16 & 11 \\
16 &  1 & 0 & 12 & 19 & 17 & 10 & 19 & 12 & 17 & 15 \\
10 &  2 & 18 & 0 & 20 & 18 & 13 & 15 & 15 & 14 & 17 \\
18 &  3 & 11 & 16 & 0 & 19 & 18 & 20 & 17 & 15 & 13 \\
18 &  4 & 16 & 11 & 20 & 0 & 20 & 20 & 12 & 10 & 14 \\
16 &  5 & 14 & 17 & 11 & 14 & 0 & 17 & 10 & 10 & 13 \\
17 &  6 & 11 & 18 & 12 & 16 & 10 & 0 & 11 & 12 & 12 \\
12 &  7 & 18 & 17 & 19 & 20 & 12 & 20 & 0 & 13 & 18 \\
13 &  8 & 10 & 16 & 18 & 16 & 14 & 16 & 13 & 0 & 14 \\
10 &  9 & 18 & 15 & 11 & 11 & 10 & 17 & 15 & 18 & 0 \\
  \end{array}
\right ]
\end{equation*}

\begin{figure}
  \centering
  \includegraphics[width=0.6\columnwidth]{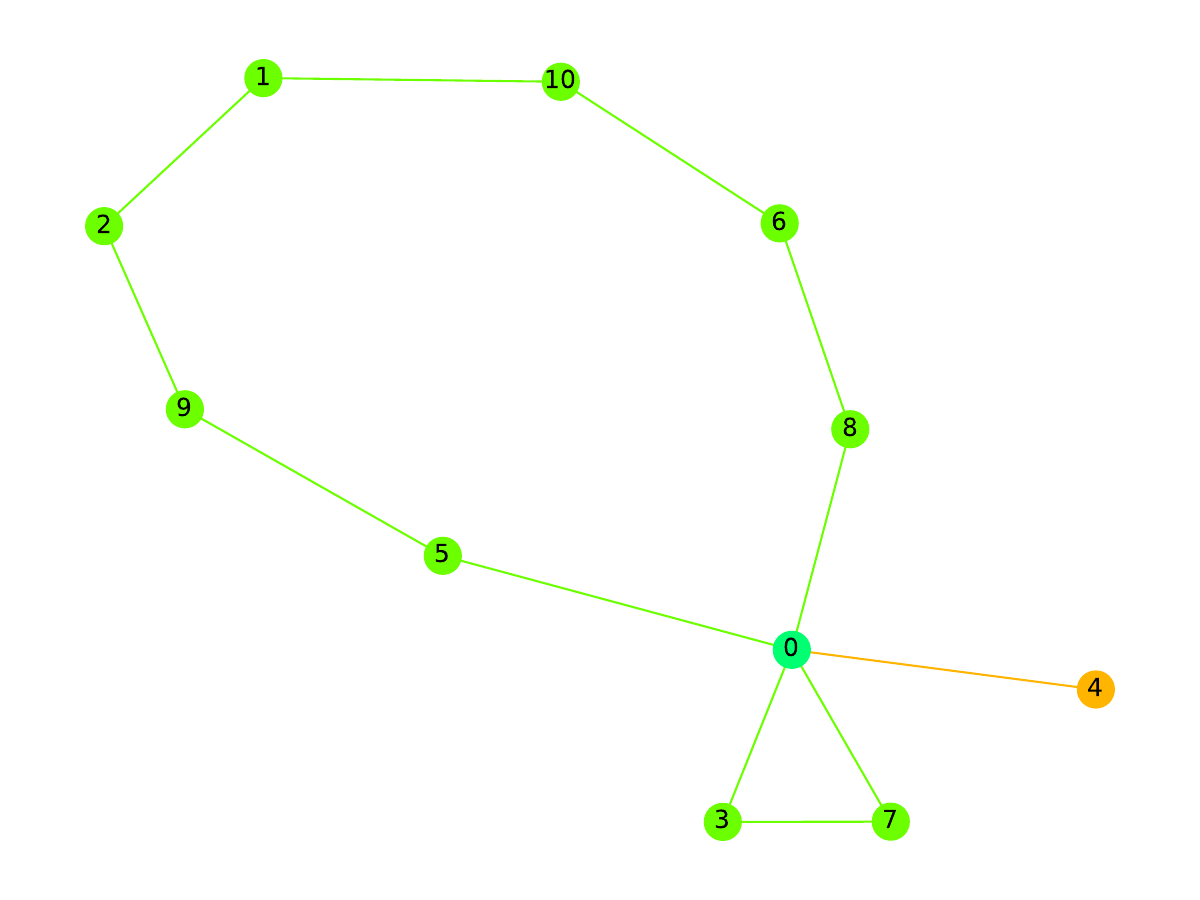}
  \caption{The route structure resulted by the example.}\label{fig:example}
\end{figure}

As can be observed from \autoref{fig:example}, there are two routes, one for each truck in use, as indicated by $|K| = 2$. Now, let's analyze the reported solution more closely:\\

\subsection{Truck 0}
For Truck 0, the optimal solution using this model suggests that only node $4$ is allocated to Truck 0, meaning $y_{4~0}=1$  and ($z_{4~0~0}= z_{4~0~1}= z_{4~0~2} =1$). \\

It's also noteworthy that $x_{400} = x_{040} = 1$. This is unusual, but it can be justified by the presence of a node with particularly high demand.\\

An interesting observation is the following: $u_{0~0~0} = 7, ~
u_{0~0~1} = 7, ~u_{0~0~2} = 6, ~u_{1~0~0} = 8, ~u_{1~0~1} = 5, ~u_{1~0~2} = 4, ~u_{2~0~0} = 4, ~u_{2~0~1} = 6, ~u_{2~0~2} = 3, ~u_{3~0~0} = 3, ~u_{3~0~1} = 2, ~u_{3~0~2} = 2, ~
u_{4~0~0} = 6, ~u_{4~0~1} = 4, ~u_{4~0~2} = 3, ~u_{5~0~0} = 6, ~u_{5~0~1} = 1, ~u_{5~0~2} = 3, ~u_{6~0~0} = 8, ~u_{6~0~1} = 7, ~u_{6~0~2} = 7, ~ u_{7~0~0} = 500, ~u_{7~0~1} = 500, ~u_{7~0~2} = 500, ~u_{8~0~1} = 4, ~u_{8~0~2} = 7, ~u_{9~0~0} = 2, ~u_{9~0~1} = 3, ~u_{9~0~2} = 2$\\

All the previously mentioned variables have non-zero values, even though only  $u_{4~0~0}, ~u_{4~0~1}$ and $u_{4~0~2}$ should be non-zero, as node $4$ is the only one allocated to truck 0. \\

\paragraph{\textbf{Issue 1}} The problem lies with constraints \eqref{eq6}, as it enforces a lower bound on variables even when they are irrelevant to the solution.

\paragraph{\textbf{Fix 1}}
To address this issue, we should replace the ill-defined inequality \eqref{eq6}, which forces every variable to take a lower bound value, even when it is irrelevant to the solution. These constraints should be replaced with the following:

\begin{align}
&q_{ip}y_{ik} \leq u_{ikp} \leq Q_p y_{ik}  & \forall i\in N_c, k\in K, p\in P\label{eq60}
\end{align}

This ensures that when node $i$ is not allocated to the route of truck $k$, the variable should not take a non-zero value for any of its compartments.

\begin{figure}
  \centering
  \includegraphics[width=0.6\columnwidth]{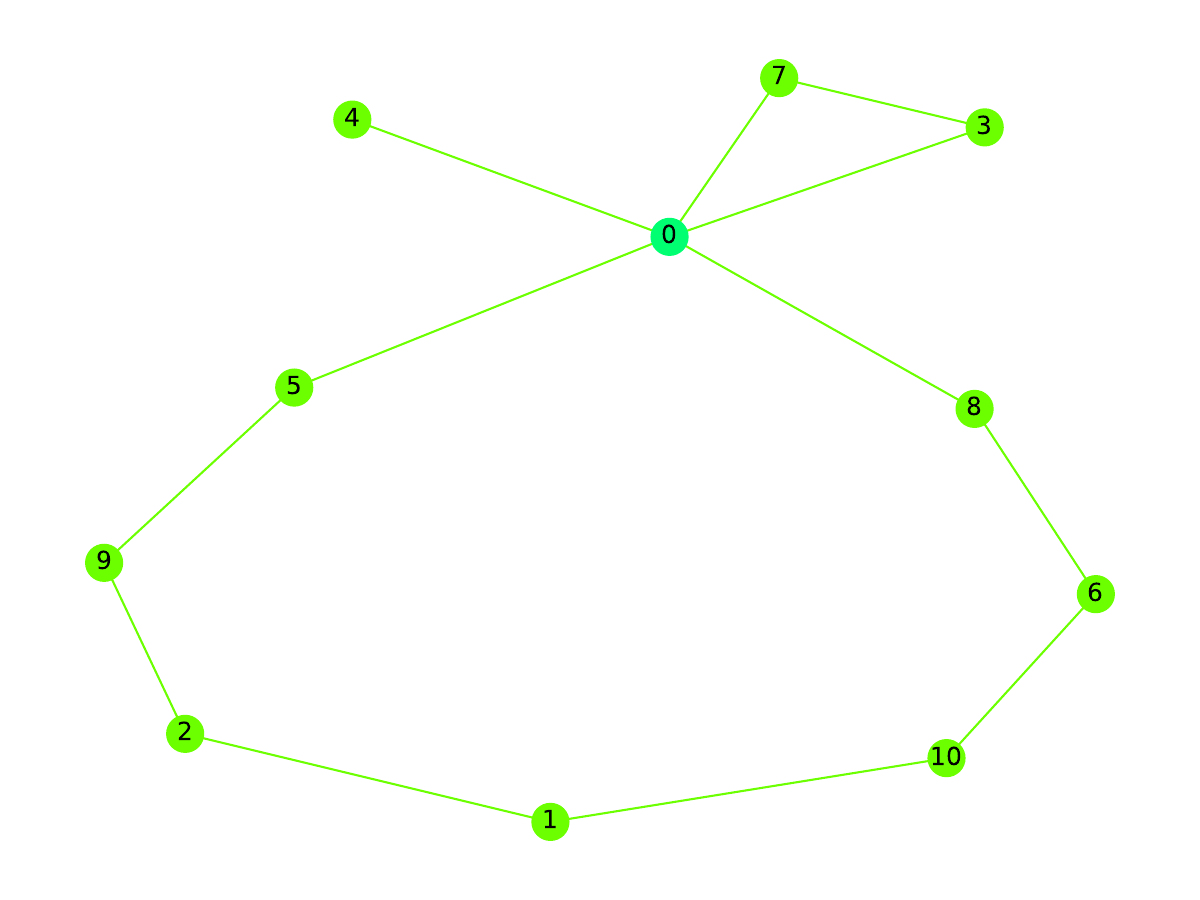}
  \caption{The new optimal solution after adding \eqref{eq60}.}\label{fig:example1}
\end{figure}

Assuming that this constraint has been added to the model and part of the issue has been resolved, the new optimal solution structure is depicted in \autoref{fig:example1}. None of the variables for truck 0 will have a non-zero value, and there are six edges incident to the depot (node 0) for truck 1.

\paragraph{\textbf{Issue 2}}
Considering constraint \eqref{eq2}, the depot is determined to be assigned to only one truck, yet six edges are incident to it.

\paragraph{\textbf{Fix 2}}

We  add the following constraint to fix it:

\begin{align}
&\sum_k y_{0k} = \sum_j \sum_k x_{0jk}  & \label{eq20}
\end{align}

\paragraph{\textbf{Issue 3}}

\begin{figure}
  \centering
  \includegraphics[width=0.6\columnwidth]{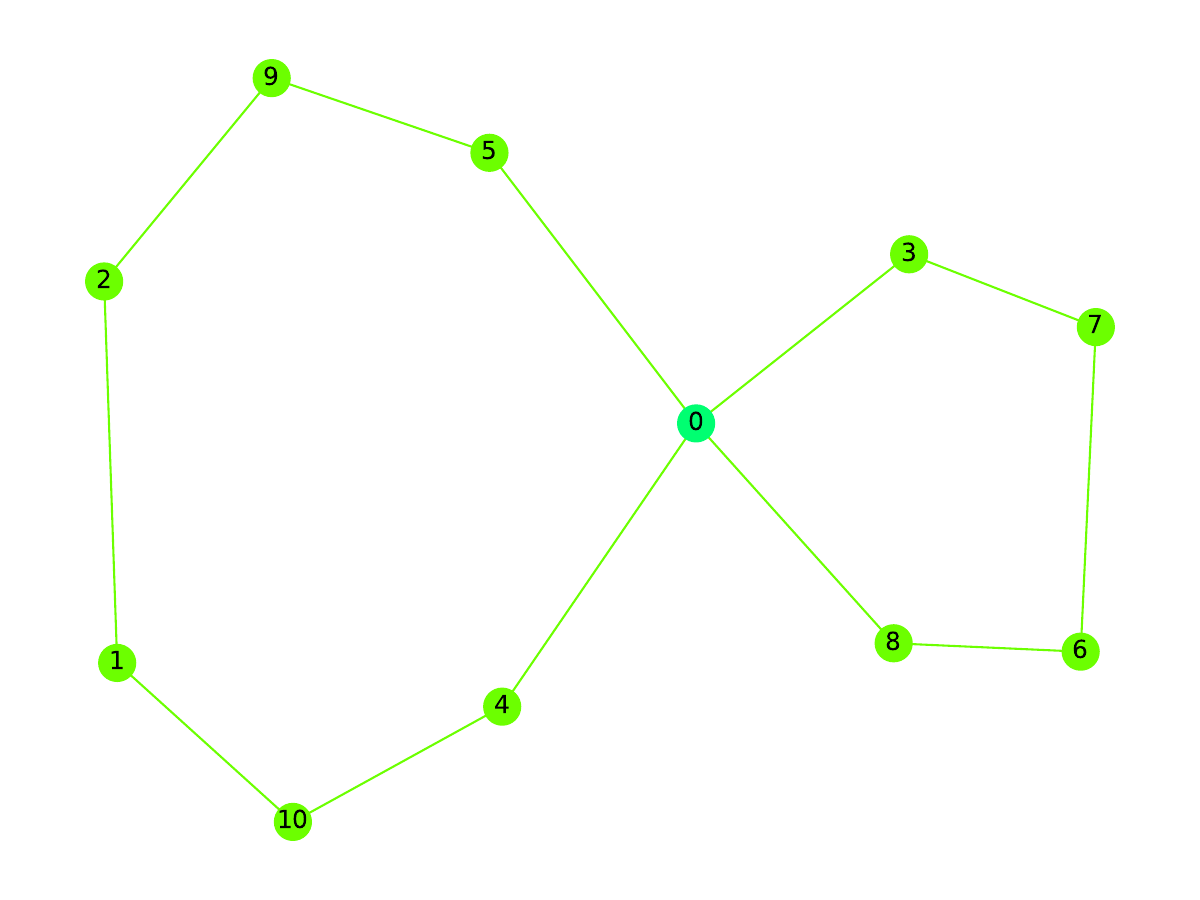}
  \caption{The new optimal solution after adding \eqref{eq20}.}\label{fig:example2}
\end{figure}

While this partially addresses the issue, it introduces a new unusual structure, as depicted in \autoref{fig:example2}. In this case, the depot (0) is randomly assigned to a second truck without any other non-zero variables for this route (0). Consequently, the left-hand side of constraint \eqref{eq20} equals 2, and two edges $(0, *)$ depart from the depot, thus satisfying constraint \eqref{eq20}. However, the resulting solution is still incorrect.

\paragraph{\textbf{Fix 3}}
As we will see in the following, a similar issue arises for the second truck (Truck 1) in the original optimal solution. The solution is to introduce a specialized form of sub-tour elimination constraints, which are exponential in number.

\subsection{Truck 1}
Now, let's revisit the original optimal solution before the addition of constraints \eqref{eq60} and \eqref{eq20}, as shown in \autoref{fig:example}. The situation in this segment (indicated by green nodes and edges) is more complex. Several issues can be highlighted in this part of the solution.

For this truck the tour is divided into two parts. The first sub-tour includes $x_{0~3~1}=x_{3~7~1}==x_{7~0~1}=1$ and the second one is comprised of $~x_{0~5~1}= ~x_{5~9~1}= ~x_{9~2~1}= ~x_{2~1~1}=x_{1~10~1}= ~x_{10~6~1}= ~x_{6~8~1}= ~x_{8~0~1}= 1$. This is a butterfly structure wherein the origin (depot) is a common point between the two sub-tours. Furthermore, we have, $z_{1~1~0}=z_{1~1~1}=z_{1~1~2}=z_{2~1~0}=z_{2~1~1}=z_{2~1~2}=z_{3~1~0}=z_{3~1~1}=z_{3~1~2}=z_{5~1~0}=z_{5~1~1}=z_{5~1~2}=z_{6~1~0}= z_{6~1~1}=z_{6~1~2}=z_{7~1~0}=z_{7~1~1}=z_{7~1~2}=z_{8~1~0}=z_{8~1~1}=z_{8~1~2}=z_{9~1~0}=z_{9~1~1}=z_{9~1~2}=z_{10~1~0}= z_{10~1~1}=z_{10~1~2}=1$ and $ y_{0~1}= y_{1~1} = y_{2~1} = y_{4~1} = y_{5~1} = y_{6~1} = y_{7~1} = y_{8~1} = y_{9~1}$.\\

\paragraph{\textbf{Issue 3'}} Let us look at the solution after the rectifications proposed so far ---i.e. \autoref{fig:example2}. The depot (node 0) is visited 4 times, which contradicts the intended purpose of this model. Normally, a return to the depot is expected to complete the tasks of a truck, which does not seem to be the case in here.

\paragraph{\textbf{Fix 3'}}
To prevent such sub-tours, an exponential number of constraints can be applied as follows. Let $S\subset V, ~ 0\in S$ and $\delta(S) = \{e=\{i.j\}| i\in S, j\in V\backslash\ S\}$
\begin{align}\label{eq:subtour}
&  \sum_{e\in \delta(S\backslash \{0\})}x_{e=\{i,j\}~k } \geq  2z_{ik}, & \forall k\in K, i \in V\backslash \ S, S\subset V, 0\in S
\end{align}

\paragraph{\textbf{Issue 4}}
The remaining variables are related to the capacity constraints, which indicate the truck load while departing a facility node. Here, $u_{0~1~0} = 21, ~u_{0~1~1} = 20, ~u_{0~1~2} = 20, ~u_{1~1~0} = 14, ~u_{1~1~1} = 13, ~u_{1~1~2} = 14, ~u_{2~1~0} = 4, ~u_{2~1~1} = 6, ~u_{2~1~2} = 3, ~u_{3~1~0} = 3, ~u_{3~1~1} = 2, ~u_{3~1~2} = 2, ~u_{4~1~0} = 6, ~u_{4~1~1} = 4, ~u_{4~1~2} = 3, ~u_{5~1~0} = 29, ~u_{5~1~1} = 24, ~u_{5~1~2} = 25, ~u_{6~1~0} = 12, ~u_{6~1~1} = 13, ~u_{6~1~2} = 10, ~u_{7~1~0} = 500, ~u_{7~1~1} = 500, ~u_{7~1~2} = 500, ~u_{8~1~0} = 6, ~u_{8~1~1} = 8, ~u_{8~1~2} = 10, ~u_{9~1~0} = 23, ~u_{9~1~1} = 23, ~ u_{9~1~2} = 22$.

The flow conservation principle is ineffective in this context. Moreover, as observed, the highest load is 29 units, which is well below the capacity limit. Evidently, the issues related to constraints \eqref{eq6} persist here as well.\\

Interestingly, some variables take values equivalent to the vehicle's compartment capacity, i.e., 500, even though given our instance the load on any compartment of a truck never comes close to this value. This discrepancy is due to the fact that constraints \eqref{eq5} are insufficient to ensure that the solver doesn't arbitrarily assign values to these variables. One might argue that this is due to the butterfly/sub-tour structure and would be resolved with the introduction of constraints \eqref{eq60}. However, this issue also affects Truck 0, even though it doesn't have a butterfly structure.

\paragraph{\textbf{Fix 4}}
This issue can only be resolved if the flow has a specific direction to follow. Even in that scenario, it would be necessary to include constraints that can manage flow conservation, more effectively.

\section{Results and Discussion}
A model is presented in \cite{OUERTANI2023102104} that serves as a reference point for measuring the quality of the adaptive genetic algorithm they propose. However, this model has several deficiencies and flaws. It can only provide a lower bound with an unknown distance from its optimal solution, making it unsuitable as a reliable reference point for assessing the quality of a heuristic or determining the optimality of solutions reported by a given heuristic. While we have proposed fixes for some issues in this model, there are still flaws for which corrections, without significant changes in the model and variables' definitions, appear challenging.

%%\subsection{Small Size Instance}
%%\subsection{Real-life Size Instances}
%\section{Conclusion}
%
%%------------------------------------------------
%\phantomsection
%\section*{Acknowledgments} % The \section*{} command stops section numbering
%
%\addcontentsline{toc}{section}{Acknowledgments} % Adds this section to the table of contents
%
%Our acknowledgement addresses the honest, respected and committed researchers all around the globe combatting against scientific exploitations, power abuse, racial/ethnic discrimination, spoofed and manipulated results, fake publications, plagiarism and all other malign  activities smearing scientific community.

%----------------------------------------------------------------------------------------
%	REFERENCE LIST
%----------------------------------------------------------------------------------------
%\section*{References}

\bibliographystyle{apalike}
\bibliography{sample}

%----------------------------------------------------------------------------------------

\end{document}